\documentclass[11pt]{article}
\usepackage{amsmath,amsfonts,amssymb,amsthm,amscd}
\usepackage{color}
\numberwithin{equation}{section} \allowdisplaybreaks
 \usepackage[all]{xy}

\newtheorem{theorem}{\sc Theorem}[section]
\newtheorem{lemma}[theorem]{\sc Lemma}
\newtheorem{proposition}[theorem]{\sc Proposition}
\newtheorem{corollary}[theorem]{\sc Corollary}
\newtheorem{definition}[theorem]{\sc Definition}
\newtheorem{example}[theorem]{\sc Example}
\newtheorem{remark}[theorem]{\sc Remark}

\newcommand{\bet}{\begin{theorem}}
\newcommand{\eet}{\end{theorem}}
\newcommand{\blm}{\begin{lemma}}
\newcommand{\elm}{\end{lemma}}
\newcommand{\bprop}{\begin{proposition}}
\newcommand{\eprop}{\end{proposition}}
\newcommand{\bcor}{\begin{corollary}}
\newcommand{\ecor}{\end{corollary}}
\newcommand{\bdf}{\begin{definition}\rm}
\newcommand{\edf}{\end{definition}}
\newcommand{\bp}{\begin{proof}}
\newcommand{\ep}{\end{proof}}
\newcommand{\bex}{\begin{example}\rm}
\newcommand{\eex}{\end{example}}
\newcommand{\bremark}{\begin{remark}\rm}
\newcommand{\eremark}{\end{remark}}

\newcommand{\norm}[1]{ \| #1 \| }

\begin{document}

\title {Index of continuous families of bounded linear operators in Banach spaces and application}

\author{ M. Berkani}

\date{}

\maketitle

\begin{abstract}

In this paper, we define an analytical index for  continuous
families of Fredholm  operators parameterized by a topological
space $\mathbb{X}$ into a Banach space $X.$ We also consider the
Weyl spectrum for continuous families of bounded linear  operators
and we study its continuity.

\end{abstract}
\footnotetext{\hspace{-7pt}2010 {\em Mathematics Subject
Classification\/}:  47A53, 47A60, 58B05
\baselineskip=18pt\newline\indent {\em Key words and phrases\/}:
   connected components, Fredholm, homotopy, index, Weyl specrtrum}

\section { Index of continuous families of Fredholm
operators} \setcounter{df}{0}

Let $X$ be an infinite dimensional Banach space, let $L(X)$  be
the Banach algebra of bounded linear operators acting on
 $X$ and let $ T \in L(X).$ We
will denote by $N(T)$ the null space of $T$, by $ \alpha(T)$ the
nullity of $T$,  by $ R(T)$ the range of $T$ and by $\beta(T) $
its defect. If the range $ R(T)$ of $T$ is closed and $ \alpha(T)
< \infty $ (resp. $\beta(T)< \infty $ ), then $T$ is called an
upper semi-Fredholm (resp. a lower semi-Fredholm) operator. A
semi-Fredholm operator is an upper or a lower semi-Fredholm
operator. In the sequel $\Phi_{+}(X)$ (resp. $\Phi_{-}(X)$) will
denote the set of upper (resp. lower) semi-Fredholm operators . If
both of $ \alpha(T) $ and $\beta(T) $ are finite, then $T $ is
called a Fredholm operator and the index of $T$  is defined by $
ind(T) = \alpha(T) - \beta(T). $

In \cite{P48}, we defined an analytical index for a continuous
family of Fredholm operators parameterized by a topological space
$\mathbb{X}$ into a Hilbert space $H,$ as a sequence of integers,
extending naturally the usual definition of the index of a single
Fredholm operator and we proved the homotopy invariance of this
index. We proved also that if $\mathbb{X}$ is a compact locally
connected space, satisfying  an homotopy condition, then the
analytical index establishes an isomorphism between the homotopy
equivalence classes of families of Fredholm operators
$[\mathbb{X}, Fred(H)]$ and the group $\mathbb{Z}^{n_c},$ where
$n_c$ is the cardinal of the connected components of $\mathbb{X},$
proving by this way a result similar to the theorem of
Atiyah-J$\ddot{a}$nich \cite[Theorem 3.40]{BB}.

The motivation of \cite{P48},  was  the construction of  an
analytical index  for continuous families of Fredholm operators
parameterized by a topological space as a sequence of integers,
 extending naturally the usual index of a single Fredholm operator,
which is an integer and avoiding the use of vector bundles.

Here, in this section, we extend the  definition of the analytical
index given in \cite{P48}  to the case of a continuous family of
Fredholm operators parameterized by a topological space
$\mathbb{X}$ into a Banach space $X.$

Define an equivalence relation $\mathcal{R}$ on the set
$fdim(X)\times fcod(X),$  where $fdim(X)$   is the set of   finite
dimensional vector subspaces of $X$  and $fcod(X)$   is the set of
finite codimension vector subspaces of $X$  by:

 $$(E_1,F_1) \mathcal{R} (E_1', F_1') \Leftrightarrow  dim E_1 -  codim F_1 = dim E_1'-codim F_1', $$

where $ dim$ stands for dimension, while $ codim $ for
codimension.

\noindent Since $X$ is an infinite dimensional vector space, then
the map: $$ \Psi: fdim(X)\times fcod(X)/\mathcal{R}  \rightarrow
\mathbb{Z},$$  defined by $\Psi( \overline{(E_1, F_1)} )=dimE_1-
codim F_1,$  where $\overline{(E_1, F_1)} $ is the equivalence
class of the couple $(E_1, F_1),$  is a bijection. Moreover $\Psi$
generate a commutative group structure on the set [$fdim(X)\times
fcod(X)]/R$ and $\psi$ is then a group isomorphism.\\

\noindent Consider now a family of Fredholm operators parametrized
by a topological space $\mathbb{X},$ that is  a continuous map $
\mathcal{T}: \mathbb{X} \rightarrow Fred(X),$  where $ Fred(X)$ is
the set of Fredholm operators, endowed with the norm topolgy of
$L(X).$ We denote by $ \mathcal{T}_x$ the image $\mathcal{T}(x)$
of an element  $ x \in \mathbb{X}.$

Define   an equivalence relation $ \textbf{\textit{c}} $ on the
space $\mathbb{X}$ by setting that $ x \,\textbf{\textit{c}}\, y,$
if and only if $ x$ and $y$ belongs to the same connected
component of $\mathbb{X}.$ Let $\mathbb{X}^{\textbf{\textit{c}}} $
be the quotient space associated to this equivalence relation, let
$ \mathcal{C}( \mathbb{X}, Fred(X))$ be the space of continuous
maps from the topological space $\mathbb{X}$ into the topological
space $ Fred(X)$ and let  the
 map: $$    q: \mathcal{C}(\mathbb{X}, Fred(X)) \rightarrow [
[fdim(X)\times fcod(X)]/\mathcal{R
}]^{\mathbb{X}^\textbf{\textit{c}}},$$

 defined by
 $q(\mathcal{T})= ( \overline{N(\mathcal{T}_{x}),R(\mathcal{T}_{x})})_{\overline{x} \in
 \mathbb{X}^{c}}$ for all  $ \mathcal{T} \in \mathcal{C}( \mathbb{X}, Fred(X)).$
Define also the map $$\Psi_ {\scriptscriptstyle\mathbb{X}}:
[[fdim(X)\times fcod(X)]/\mathcal{R}
]^{\mathbb{X}^{\textbf{\textit{c}}}}\rightarrow
\mathbb{Z}^{n_\textit{\textbf{c}}}$$
 by setting $\Psi_ {\scriptscriptstyle\mathbb{X}}((Y_{\overline{x}}) _{
\overline{x} \in \mathbb{X}^{\textbf{\textit{c}}} })= (\Psi
(Y_{\overline{x}})) _{\overline{x} \in
\mathbb{X}^{\textbf{\textit{c}}}} $ for all $ (Y_{\overline{x}})
_{ \overline{x} \in \mathbb{X}^{\textbf{\textit{c}}} } \in
[fdim(X)\times fcod(X)]/\mathcal{R}
]^{\mathbb{X}^{\textbf{\textit{c}}}}.$ Here
$n_{\textbf{\textit{c}}} $ stands for the cardinal of the
connected components of the topological space $\mathbb{X},$
assuming that the space  $\mathbb{X}$ has at most a countable set
of  connected components.

 \bdf \label{index}  The analytical  index (or simply the index) of  a family of Fredholm operators $ \mathcal{T} : \mathbb{X} \rightarrow Fred(X),$
   parameterized by a topological space $\mathbb{X}$    is defined by $ ind(\mathcal{T})=
\psi_{\scriptscriptstyle  \mathbb{X}}( q(\mathcal{T})).$
 \edf

Explicitly, we have:   $ ind(\mathcal{T})=
\psi_{\scriptscriptstyle \mathbb{\mathcal{T}}} ((
\overline{N(\mathcal{T}_{x}),R(\mathcal{T}_{x})})_{\overline{x}
\in
 \mathbb{X}^{c}}) = ( ind(\mathcal{T}_x))_{\overline{x}
\in \mathbb{X}^{c}}. $

\noindent Thus the index of  a family of Fredholm operators $
\mathcal{T}$ is a sequence of integers in
$\mathbb{Z}^{n_\textit{\textbf{c}}},$ which may be a finite
sequence or infinite sequence, depending on the cardinal of the
connected components of $\mathbb{X}.$

 \bet The index of a continuous family of  Fredholm  operators $\mathcal{T}$  parameterized
by a topological space $\mathbb{X} $ is well defined as an element
of $ \mathbb{Z}^{ n_c}.$ In particular if $\mathbb{X} $ is reduced
to a single element, then the  index of $\mathcal{T}$  is equal to
the usual index of the Fredholm operator $\mathcal{T}.$ \eet

\bp From the usual properties of the index \cite[Theorem
5.17]{KA}, we know that two Fredholm operators located in the same
connected component of the set of Fredholm operators have the same
index. Moreover, as $ \mathcal{T} $ is continuous, the image of a
connected component of the topological space $ \mathbb{X},$ is
included in a connected component of the set of Fredholm
operators.  This shows that the index of a family of Fredholm
operators  is well defined, and it is clear that if $ \mathbb{X}$
is reduced to a single element, the  index of $\mathcal{T}$
defined  here is equal to the usual index of the single Fredholm
operator $T.$

\ep

\bdf A continuous  family $\mathcal{K}$  from  $\mathbb{X}$ to
$L(X)$ is said to be compact if $\mathcal{K}_x$ is compact for all
$ x \in X.$ The family $\mathcal{K}$  is said to be of finite rank
  if  $\mathcal{K}_x$ is of finite rank for all $ x \in X.$

\edf

\bprop i) Let $\mathcal{T} \in  \mathcal{C}( \mathbb{X}, Fred(X))
$ and  let $ \mathcal{K}$ be a  continuous  compact   family from
$\mathbb{X}$ to $L(X).$ Then $ \mathcal{T}+\mathcal{K} $ is a
Fredholm family and $ ind(\mathcal{T})=
ind(\mathcal{T}+\mathcal{K}).$

ii) Let $\mathcal{S,T } \in  \mathcal{C}( \mathbb{X}, Fred(X))$ be
two Fredholm families, then the family $\mathcal{ST}$ defined by
$(\mathcal{ST})_x= \mathcal{S}_x\mathcal{T}_x$ is a Fredholm
family and $ ind(\mathcal{ST})= ind(\mathcal{S})+
ind(\mathcal{T}).$

\eprop

\bp This is clear from the usual properties of Fredholm operators.

\ep

\bet \label{Ideal} Assume that $\mathbb{X}$ is a compact
topological space. Then the  set $ \mathcal{K} \mathcal{C}(
\mathbb{X}, L(X))$ of continuous compact families from
$\mathbb{X}$ to $L(X)$ is a closed ideal in the Banach algebra
$\mathcal{C}( \mathbb{X}, L(X))$  \eet

\bp  Recall that $\mathcal{C}( \mathbb{X}, L(X))$ is a unital
algebra with the usual properties of addition, scalar
multiplication and multiplication defined by:

$ (\lambda \mathcal{S+ T})_x= \lambda \mathcal{S}_x +
\mathcal{T}_x,   (\mathcal{ST})_x=
\mathcal{S}_x\mathcal{\mathcal{T}}_x,
 \forall \mathcal{(S, T)} \in \mathcal{C}( \mathbb{X}, L(X))^2,  \forall
x \in \mathbb{X},  \forall  \lambda \in \mathbb{C}.$

The unit element of $\mathcal{C}( \mathbb{X}, L(X))$ is the
constant function $\mathcal{I}$ defined by $ \mathcal{I}_x= I $
the identity of $X,$  for all $ x \in \mathbb{X}.$ Moreover as
$\mathbb{X}$ is compact, then if we set $ \norm{\mathcal{T}}=
sup_{x\in \mathbb{X}} \norm {\mathcal{T}_x}, \forall \mathcal{T}
\in \mathcal{C}( \mathbb{X}, L(X)),$  then $\mathcal{C}(
\mathbb{X}, L(X))$ equipped with this norm is a Banach algebra.
Similarly $\mathcal{C}( \mathbb{X}, L(X)/K(X))$ equipped with the
norm  $ \norm{\mathcal{T}}= sup_{x\in \mathbb{X}} \norm
{P\mathcal{T}_x}$ is a unital Banach algebra, where  $ P: L(X)
\rightarrow L(X)/K(X)$ is  the usual projection from $ L(X)$ onto
the Calkin algebra $L(X)/K(X).$

 It is clear that $ \mathcal{K} \mathcal{C}( \mathbb{X}, L(X))$ is
 an ideal of  $  \mathcal{C}( \mathbb{X}, L(X)).$ Assume now that
 $ (\mathcal{T}_n)_n$ is a sequence in $ \mathcal{K}\mathcal{C}( \mathbb{X}, L(X))$
 converging  in $ \mathcal{C}( \mathbb{X}, L(X)$
 to $\mathcal{T}.$ Then  $ ((\mathcal{T}_n)_x)_n$ converges to $\mathcal{T}_x,$ as each $(\mathcal{T}_n)_x
 $ is compact, then $\mathcal{T} \in \mathcal{K}\mathcal{C}( \mathbb{X}, L(X)).$
\ep

\bremark  In the same way as in the case of the Calkin algebra,
Theorem \ref{Ideal} generates the new Banach algebra $\mathcal{C}(
\mathbb{X}, L(X))/ \mathcal{K} \mathcal{C}( \mathbb{X}, L(X)).$
Moreover, there is a natural injection $\overline{\Pi}:
\mathcal{C}( \mathbb{X}, L(X))/ \mathcal{K} \mathcal{C}(
\mathbb{X}, L(X)) \rightarrow  \mathcal{C}( \mathbb{X},
L(X)/K(X))$ defined by $\overline{\Pi}( \overline{\mathcal{T}})=
P\mathcal{T},$ where $ \overline{\mathcal{T}}$ is the equivalence
class of the element $\mathcal{T}$ of $ \mathcal{C}( \mathbb{X},
L(X))$ in $\mathcal{C}( \mathbb{X}, L(X))/ \mathcal{K}
\mathcal{C}( \mathbb{X}, L(X))$ and $ P: L(X) \rightarrow
L(X)/K(X)$ is the natural projection.

\vspace{3mm} \noindent \textbf{\underline{Open questions:}}

\begin{enumerate}

\item Given an element $ \mathcal{S} \in \mathcal{C}(\mathbb{X},
L(X)/K(X)),$  does there exist a continuous family $ \mathcal{T}
\in \mathcal{C}(\mathbb{X}, L(X))$ such that $ \overline{\Pi}(
\overline{\mathcal{T}})= \mathcal{S}?$

\item Let $H$ be a Hilbert space and  $ \mathcal{K}$ be a
continuous compact family in $ \mathcal{C}(\mathbb{X},L(H).$ Does
there exists a sequence   $ (\mathcal{K}_n)_n$ of  continuous
finite rank families in $ (\mathcal{C}(\mathbb{X},L(H))$ such that
$ \lim \limits _{n\rightarrow \infty }{}\mathcal{K}_n =
\mathcal{K}?$

\end{enumerate}

\eremark

\bet
 Assume that $\mathbb{X}$ is a compact topological space  and  let $\mathcal{T} \in \mathcal{C}(\mathbb{X}, L(X)).$ Then $ \mathcal{T}$ is a
 Fredholm family if and only if $P\mathcal{T}$ is  invertible in the Banach algebra
 $ \mathcal{C}(\mathbb{X}, L(X)/K(X)).$

\eet

\bp Assume that $\mathcal{T}$ is a Fredholm family, then for all
$x \in X, \mathcal{T}_x$  is a Fredholm operator. Thus
$P\mathcal{T}_x$ is invertible in $L(X)/K(X).$ Let
$(P\mathcal{T}_x)^{-1}$ be its inverse, then the family
$(P\mathcal{T})^{-1}$ defined by  $ (P\mathcal{T})^{-1}(x)=
(P\mathcal{T}_x)^{-1}$ is a continuous family, because the
inversion is a continuous map in the Banach algebra $L(X)/K(X),$
and $(P\mathcal{T})^{-1}$ is the inverse of $P\mathcal{T}$ in the
Banach algebra $ \mathcal{C}(\mathbb{X}, L(X)/K(X)).$

Conversely if $ P\mathcal{T}$  is  invertible in the Banach
algebra
 $ \mathcal{C}(\mathbb{X}, L(X)/K(X)),$ then there exists $S \in   \mathcal{C}(\mathbb{X}, L(X)/K(X))$
such that $ (P\mathcal{T})\mathcal{S}= \mathcal{S}(P\mathcal{T})=
\overline{\mathcal{I}},$  where $\overline{\mathcal{I}}$ is
defined by $\overline{\mathcal{I}}_x= PI,$ for all $ x \in
\mathbb{X}.$ Thus $(P\mathcal{T}_x) \mathcal{S}_x=
\mathcal{S}_x(P\mathcal{T}_x)= PI.$ Thus $P\mathcal{T}_x $ is
invertible in the Calkin algebra $L(X)/K(X),$
 $\mathcal{T}_x$ is a Fredholm operator and $\mathcal{T} \in \mathcal{C}(\mathbb{X},
Fred(X)).$

\ep


\bdf Let $ \mathcal{S,T} $   be in $\in \mathcal{C}( \mathbb{X},
Fred(X)). $ We will say that $\mathcal{S}$ and $\mathcal{T}$ are
Fredholm homotopic, if there exists a map $\Phi:  [0,1]\times X
\rightarrow L(X)$ such that $ \Phi(0,x)= \mathcal{S}_x, \,
\Phi(1,x)= \mathcal{T}_x $ and $\Phi(t,x)$ is a Fredholm operator,
for all $ (t,x) \in [0,1]\times \mathbb{X}.$ \edf

\bet Let $\mathcal{ S, T}$ be two Fredholm homotopic elements of
\,$\mathcal{C}( \mathbb{X}, Fred(X)). $  Then  $ind(\mathcal{T})=
ind(\mathcal{S}).$ \eet

\bp Since $\mathcal{S}$ and $\mathcal{T}$ are Fredholm homotopic,
there exists a continuous map $ h: \mathbb{X} \times [0, 1]
\rightarrow Fred(X) $  such that  $h(x,0)=\mathcal{S}(x)$ and
$h(x,1)=T(x)$ for all $x\in\mathbb{X}.$ For a fixed $x \in
\mathbb{X} ,$ the map $ h_x: [0,1] \rightarrow Fred(X),$ defined
by $h_x(t)= h(x,t)$ is a continuous path
 in $ Fred(X)$ linking  $\mathcal{S}_x$ to $\mathcal{T}_x.$  Since $[0,1]$ is connected and since the
 usual index is constant on connected sets of $Fred(X)$, then $ind(\mathcal{S}_x)=
 ind(\mathcal{T}_x).$ So  $q(\mathcal{S})=q(\mathcal{T})$ and then  $ind(\mathcal{S})= ind(\mathcal{T}).$

\ep

\bet \label {index-continuity} Let $\mathbb{X}$ be a compact
topological  space. Then the index is a continuous locally
constant function  from $ \mathcal{C}( \mathbb{X}, Fred(X)) $ into
the group $\mathbb{Z}^{ n_\textit{\textbf{c}}}.$

\eet

\bp Let $\mathcal{T} \in \mathcal{C}( \mathbb{X}, Fred(X)),$  then
$\forall x \in  \mathbb{X}, \exists \, \epsilon_x > 0, $ such that
$B(\mathcal{T}_x, \epsilon_x) \subset Fred(X),$  because $Fred(X)$
is open in $L(X).$ Then the index is constant on $B(\mathcal{T}_x,
\epsilon_x),$ because $B(\mathcal{T}_x, \epsilon_x)$ is connected.
We have $ \mathbb{X} \subset \bigcup_{x\in
\mathbb{X}}\mathcal{T}^{-1}( B(\mathcal{T}_x,
\frac{\epsilon_{x}}{2})).$ Since $\mathbb{X}$ is compact, there
exists $ x_1,...x_n$ in $\mathbb{X}$ such that $\mathbb{X} \subset
\bigcup_{i=1}^n\mathcal{ T}^{-1}(B(T_{x_i},
\frac{\epsilon_{x_i}}{2})).$ Let $ \epsilon= min \{
\frac{\epsilon_{x_i}}{2} | 1\leq i \leq n \} $ the minimum of the
$ \frac{\epsilon_{x_i}}{2}, 1\leq i \leq n,$ and  let $S \in
\mathcal{C}( \mathbb{X}, Fred(X)) $ such that $ || \mathcal{T-S}||
< \frac{\epsilon}{2}.$ If $ x \in \mathbb{X},$ then $ ||
\mathcal{T}_x-\mathcal{S}_x|| < \frac{\epsilon}{2}$ and there
exists $ i, 1\leq i \leq n,$ such that $ x \in\mathcal{
T}^{-1}(B(T_{x_i}, \frac{\epsilon_{x_i}}{2})).$ Then $ ||
\mathcal{S}_x -\mathcal{T}_{x_i}|| \leq ||
\mathcal{S}_x-\mathcal{T}_x|| + ||\mathcal{T}_x-
\mathcal{T}_{x_i}|| < \epsilon/2 + \epsilon_{x_i}/2 \leq
\epsilon_{x_i}.$ So $ind(\mathcal{S}_x)= ind(\mathcal{T}_{x_i})=
ind(\mathcal{T}_x).$ Hence the index is a locally constant
function, in particular it is a continuous function. \ep

\bet  Let $\mathbb{X}$ be a compact topological space.  Then the
set $ \mathcal{C}( \mathbb{X}, Fred(X))$ is an open subset of the
Banach algebra $ \mathcal{C}( \mathbb{X}, L(X))$ endowed with the
uniform norm $ \norm{ \mathcal{T}}= sup_{x\in  \mathbb{X}}
\norm{\mathcal{T}_x}.$ \eet

\bp Let $\mathcal{T} \in \mathcal{C}( \mathbb{X}, Fred(X)),$  then
$\forall x \in \mathbb{X}, \exists \, \epsilon_x
> 0, $ such that $B(\mathcal{T}_x, \epsilon_x) \subset Fred(X).$   We have $ \mathbb{X} \subset
\bigcup_{x\in \mathbb{X}} \mathcal{T}^{-1}( B(\mathcal{T}_x,
\frac{\epsilon_{x_i}}{2})).$ Since $\mathbb{X}$ is compact, there
exists $ x_1,...x_n$ in $\mathbb{X}$ such that $\mathbb{X} \subset
\bigcup_{i=1}^n T^{-1}(B(\mathcal{T}_{x_i},
\frac{\epsilon_{x_i}}{2})).$ Let $ \epsilon= min \{
\frac{\epsilon_{x_i}}{2} | 1\leq i \leq n \} $ the minimum of the
$ \frac{\epsilon_{x_i}}{2}, 1\leq i \leq n.$ Let $\mathcal{S} \in
\mathcal{C}( \mathbb{X}, L(X)) $ such that $ || \mathcal{T-S}|| <
\frac{\epsilon}{2}.$ If $ x \in \mathbb{X},$ then $ ||\mathcal{
T}_x-\mathcal{S}_x|| < \frac{\epsilon}{2}$ and there exists $ i,
1\leq i \leq n,$ such that $ x \in \mathcal{T}^{-1}(B(T_{x_i},
\epsilon_{x_i})).$ Then  $ || \mathcal{S}_x -\mathcal{T}_{x_i}||
\leq || \mathcal{S}_x-\mathcal{T}_x||  + ||\mathcal{T}_x-
\mathcal{T}_{x_i}|| < \epsilon/2 + \epsilon_{x_i}/2 \leq
\epsilon_{x_i}$ and $\mathcal{S}_x$ is a Fredholm operator.
Therefore $\mathcal{S} \in \mathcal{C}( \mathbb{X}, Fred(X))$ and
$\mathcal{C}( \mathbb{X}, Fred(X))$ is open in $ \mathcal{C}(
\mathbb{X}, L(X)).$

Alternatively, we can see that $ \mathcal{C}( \mathbb{X},
Fred(X))= \Pi^{-1}( (\mathcal{C}( \mathbb{X}, L(X)/K(X)) ^{inv}),$
where $(\mathcal{C}( \mathbb{X}, L(X)/K(X)) ^{inv}$ is the open
group of invertible elements of the unital Banach algebra
$\mathcal{C}( \mathbb{X}, L(X)/K(X))$ and $ \Pi: \mathcal{C}(
\mathbb{X}, L(X))\rightarrow \mathcal{C} ( \mathbb{X}, L(X)/K(X))$
is the map defined by $ \Pi(\mathcal{T})= P\mathcal{T},$ for all $
\mathcal{T} \in \mathcal{C}( \mathbb{X}, L(X)/K(X)).$

\ep

\bcor \label{coropen}  Let $\mathbb{X}$ be a compact topological
space and let $ p \in \mathbb{Z}^{n_c}.$ Then the set $
\mathcal{C}_p( \mathbb{X}, Fred(X))$ of the continuous Fredholm
families of index $p$ is an open subset of the Banach algebra $
\mathcal{C}( \mathbb{X}, L(X)).$

\ecor

\bp As the set $ \{p\} $ is open in $\mathbb{Z}^{n_c}$ and the
index is a continuous function,  from Theorem
\ref{index-continuity}, it follows that $ \mathcal{C}_p(
\mathbb{X}, Fred(X))$ is an open subset of the Banach algebra $
\mathcal{C}( \mathbb{X}, L(X)).$

\ep

\bet

Let  $\mathcal{\mathcal{T}} \in \mathcal{C}( \mathbb{X}, L(X))$
and $f$ an analytic function   in a
 neighborhood  of the  spectrum $\sigma(\mathcal{T})$  of   $\mathcal{T}$ which
 is non-constant on any connected component of
 $\sigma(\mathcal{T})$. Then:

 \begin{enumerate}

\item  $f(\sigma_{F}(\mathcal{T}))=\sigma_{F}(f(\mathcal{T})), $
in particular $f(\mathcal{T})$ is a Fredholm family if and only if
$f(\lambda) \neq 0,$ for all $ \lambda \in
\sigma_{F}(\mathcal{T}).$

\item  If $f(\mathcal{T})$ is a Fredholm family,  then $
ind(f(\mathcal{T}))= (  \sum \limits _{n \in \mathbb{N}}
n\alpha_{n,x} )_{\overline{x} \in \mathbb{X}^{c}},$ where for all
$ x \in \mathbb{X}, \alpha_{n,x}$ is the number of zeros of $f$ on
the set $ \{\lambda \in \sigma(\mathcal{T}_x) \mid \mathcal{T}_x
-\lambda I \in Fred(X) \, and \, ind(\mathcal{T}_x -\lambda I)= n
\},$ counted according to their multiplicity.

\end{enumerate}

\eet

\bp  \begin{enumerate} \item As we have $\sigma_{F}(\mathcal{T})=
\bigcup \limits _{x \in \mathbb{X}}{} \sigma_F(\mathcal{T}_x),$
and $f$ is an analytic function in a neighborhood  of the
spectrum $\sigma(\mathcal{T})= \bigcup \limits _{x \in
\mathbb{X}}{} \sigma(\mathcal{T}_x),$ then $f$ is analytic in a
neighborhood $\sigma(\mathcal{T}_x).$ From \cite [Theorem 1]{GL},
we obtain
$f(\sigma_{F}(\mathcal{T}_x))=\sigma_{F}(f(\mathcal{T}_x)).$ Thus:
$$f(\sigma_{F}(\mathcal{T}))=f(\bigcup \limits _{x \in
\mathbb{X}}{} \sigma_F(\mathcal{T}_x))= \bigcup \limits _{x \in
\mathbb{X}}{} f(\sigma_F(\mathcal{T}_x)))=  \bigcup \limits _{x
\in \mathbb{X}}{} \sigma_F(f(\mathcal{T}_x)))=
\sigma_{F}(f(\mathcal{T})).$$

\item  As $f$ is non-constant on any connected component of
$\sigma(\mathcal{T}), $ then it has a finite number of zeros in
$\sigma(\mathcal{T}),$ say $\lambda_1,... \lambda_k,$ with
multiplicities $ n_1,n_2,...,n_k. $ This imply in particular that
the series    $ \sum \limits _{n \in \mathbb{N}} n \alpha_{n,x} $
are finite sums, for all $ x \in \mathbb{X}.$

 Moreover, there exists a non-vanishing analytic function $g$
 on $ \sigma(\mathcal{T}),$  such that $f(\lambda)= g(\lambda) \prod \limits _{i=1}^{n} (\lambda -
\lambda_i)^{n_i}.$   Thus $g(\mathcal{T})$ is invertible, $
ind(g(\mathcal{T}))=0$ and   $ ind(f(\mathcal{T}))= ind(\prod
\limits _{i=1}^{n} (\mathcal{T} - \lambda_i \mathcal{I})^{n_i})=
\sum \limits _{i=1}^{n} n_i \,ind(\mathcal{T} - \lambda_i
\mathcal{I}) = (\sum \limits _{i=1}^{n}  n_i \, ind(\mathcal{T}_x-
\lambda_i I))_{\overline{x} \in \mathbb{X}^{c}}= (\sum \limits
_{i=1}^{n} ind(g(\mathcal{T}_x)(\mathcal{T}_x- \lambda_i
I)^{n_i}))_{\overline{x} \in \mathbb{X}^{c}}=
(ind(f(\mathcal{T}_x))_{\overline{x} \in \mathbb{X}^{c}}, $
because $g(\mathcal{T}_x)$ is invertible,   for all $ x \in
\mathbb{X}.$  Applying \cite[Theorem 1, c]{GL} to
$f(\mathcal{T}_x)$, we obtain $ind(f(\mathcal{T}_x))= \sum \limits
_{n \in \mathbb{N}} n\alpha_{n,x}.$ Observe that by the
Fredholmness of $f(\mathcal{T}),$   the integer  $ \sum \limits
_{n \in \mathbb{N}} n \alpha_{n,x} $ is constant on
$\overline{x},$ for all $  x \in \mathbb{X}.$ Thus we get the
desired result.

\end{enumerate}

\ep

\section{ On the continuity of the  Weyl spectrum}

The definition of an index for continuous families of Fredholm
operators as a sequence of integers, enable us to define the Weyl
spectrum  $\sigma_{W}(\mathcal{T})$ for a continuous family of
bounded linear operators $ \mathcal{T}.$ In \cite{OB}, upper and
lower semi-continuity properties of the map $T \rightarrow
\sigma_{W}(T),$ were studied in the case of single bounded linear
operators. Here, we study similar questions for continuous
families of bounded linear operators. For the definitions of upper
semi- continuity, lower semi-continuity, continuity , lower limit
($\liminf$), upper limit ($\limsup$) and limit ($\lim$) of sets,
we refer the reader to \cite {KU}.

\bdf  Let $\mathcal{T} \in \mathcal{C}( \mathbb{X}, L(X)). $ Then
$\mathcal{T}$ is called a Weyl   family
 if it is a
Fredholm   family of index $0$. \\
  \noindent The  Weyl
spectrum $\sigma_{W}(\mathcal{T})$ of $\mathcal{T}$ is defined by
$\sigma_{W}(\mathcal{T})= \{ \lambda \in \mathbb{C}: \mathcal{T}-
\lambda \mathcal{I} $ is not a Weyl family \}.

\edf

It's easily see that  $\sigma_{W}(\mathcal{T})=  \bigcup \limits
_{x \in \mathbb{X}}{} \sigma_{W}(\mathcal{T}_x),$  see
\cite[Theore 3.5]{P48}  for more details.

\bet\label{uppersemi} The map $\mathcal{T} \rightarrow
\sigma_{W}(\mathcal{T}), $ is an upper semi-continuous function
from $\mathcal{C}( \mathbb{X}, L(X)) $ into closed subsets of
$\mathbb{C}.$

\eet

\bp Since the index is a continuous function, then for all
$\mathcal{T}  \in \mathcal{C}( \mathbb{X}, L(X)),
\sigma_{W}(\mathcal{T}) $  is closed.

Now, let $ \mathcal{T} \in  \mathcal{C}( \mathbb{X},L(X))$ and let
$ \mathcal{T}_n$ in $\mathcal{C}( \mathbb{X}, L(X)) $ such that $
\lim \limits _{n\rightarrow \infty} \mathcal{T}_n = \mathcal {T}.$
If $ \lambda \notin \sigma_{W}(\mathcal{T}), $ then $\mathcal{T} -
\lambda I $ is a Weyl family. From Corollary \ref{coropen}, there
exists $ \eta > 0, $ such that if $ \mathcal{S} \in \mathcal{C}(
\mathbb{X}, L(X)), $ satisfy $\norm { \mathcal{T} - \lambda I -
\mathcal{S}} < \eta, $  then $ \mathcal{S} $ is a Weyl family.
Since $ \lim \limits _{n\rightarrow \infty} \mathcal{T}_n =
\mathcal {T},$ there exists an integer $N$ such that for all
integer $n \geq N,$   we have $\norm { \mathcal{T}_n - \lambda I -
(\mathcal{T}- \lambda I)} < \frac{\eta}{2}.$ If $ \mid \lambda  -
\mu \mid < \frac{\eta}{2},$ then  $\norm { \mathcal{T} - \lambda I
- (\mathcal{T}_n- \mu I)} < \eta.$ Hence $\mathcal{T}_n- \mu I$ is
a Weyl family. This imply that $ \lambda \notin
\limsup\limits_{n\rightarrow \infty} \sigma_W( \mathcal{T}_n)$ and
so $\limsup\limits_{n\rightarrow \infty} \sigma_W( \mathcal{T}_n)
\subset \sigma_W(\mathcal{T}).$ This proves that the map
$\mathcal{T} \rightarrow \sigma_{W}(\mathcal{T}) $ is  upper
semi-continuous.

\ep

\bet \label{limspectrum} Let $ \mathcal{T}_n , \mathcal{T}$  in
$\mathcal{C}( \mathbb{X}, L(X)) $  such that $ \lim \limits
_{n\rightarrow \infty} \mathcal{T}_n = \mathcal {T}.$ If $ \lim
\limits _{n\rightarrow \infty}  \sigma(\Pi(\mathcal{T}_n))=
\sigma(\Pi(\mathcal {T})),$ then $ \lim \limits _{n\rightarrow
\infty} \sigma_{W}(\mathcal{T}_n)= \sigma_{W}(\mathcal{T}).$

\eet

\bp Recall that  $ \Pi: \mathcal{C}( \mathbb{X}, L(X))\rightarrow
\mathcal{C} ( \mathbb{X}, L(X)/K(X))$ is the map defined by $
\Pi(\mathcal{T})= P\mathcal{T},$ for all $ \mathcal{T} \in
\mathcal{C}( \mathbb{X}, L(X)/K(X)),$  where $ P: L(X) \rightarrow
L(X)/K(X)$ is the natural projection.

 Since the map $\mathcal{T} \rightarrow
\sigma_{W}(\mathcal{T}) $ is an upper semi-continuous function, it
is enough to prove that $\sigma_W(\mathcal{T}) \subset \liminf
\limits _{n\rightarrow \infty} \sigma_{W}(\mathcal{T}_n).$ Suppose
that $ \lambda \notin \liminf \limits _{n\rightarrow \infty}
\sigma_{W}(\mathcal{T}_n),$ so that there is a neighborhood $V$ of
$ \lambda $ that does not intersect infinitely many
$\sigma_W(\mathcal{T}_n).$ Since $ \sigma( \Pi(T_n)) \subset
\sigma_W(\mathcal{T}_n),$ then $ V$ does not intersect infinitely
many  $\sigma(\Pi(\mathcal{T}_n)).$ Hence $ \lambda \notin  \lim
\limits _{n\rightarrow \infty} \sigma(\Pi(\mathcal{T}_n))=
\sigma(\Pi(\mathcal {T}))$  and then $ \mathcal{T}- \lambda
\mathcal{I} $ is a Fredholm family. As the index is a continuous
function, we have  $ind(\mathcal{T}- \lambda \mathcal{I})= 0$ and
$\lambda  \notin \sigma_W(\mathcal{T}).$

\ep

\bet Let $ \mathcal{T}_n , \mathcal{T}$  in \,  $\mathcal{C}(
\mathbb{X}, L(X)) $ such that $ \lim \limits _{n\rightarrow
\infty} \mathcal{T}_n = \mathcal {T}.$ Then $ \lim \limits
_{n\rightarrow \infty} \sigma_{W}(\mathcal{T}_n)=
\sigma_{W}(\mathcal{T}),$ in each of the following cases:
\begin{enumerate}

\item  $\mathcal{T}_n \mathcal{T} = \mathcal{T} \mathcal{T}_n,$
for each $n.$

\item $\sigma(\mathcal{T})$  is totally disconnected.

\item  $X$ is Hilbert space  and  $ \forall \, n,
\mathcal{T}_{n,x}$ and $\mathcal{T}_x$ are normal operators for
all $x \in \mathbb{X}.$

\end{enumerate}

\eet

\bp As  $\sigma_W(\mathcal{T})= \bigcup \limits _{x \in
\mathbb{X}}{} \sigma_W(\mathcal{T}_x) $ and
$\sigma((\Pi(\mathcal{T}))= \bigcup \limits _{x \in \mathbb{X}}{}
\sigma(\Pi(\mathcal{T}_x)),$  one of the following conditions
holds:

\begin{enumerate}

\item  $\mathcal{T}_{n,x} \mathcal{T}_x = \mathcal{T}_x
\mathcal{T}_{n,x},$ for each $n$  and for all $ x \in \mathbb{X},$

\item $\sigma(\mathcal{T}_x)$  is totally disconnected,  for all
$ x \in \mathbb{X}.$

\item  $X$ is Hilbert space  and  $ \forall \, n $ and for all   $
x \in \mathbb{X},   \mathcal{T}_{n,x}$ and $\mathcal{T}_x$ are
normal operators.

\end{enumerate}

From \cite[Corollary, p. 209]{OB}, in each of this cases, we have:
$ \lim \limits _{n\rightarrow \infty}
\sigma(\Pi(\mathcal{T}_{n,x}))= \sigma(\Pi(\mathcal {T}_x)), $ for
all  $x \in \mathbb{X}.$ From \cite [Theorem2]{OB}, we have  $
\lim \limits _{n\rightarrow \infty} \sigma_w(\mathcal{T}_{n,x})=
\sigma_w(\mathcal {T}_x)$ and from \cite[ Formula 3a, p.336]{KU},
we have  $\bigcup \limits _{x \in \mathbb{X}}{}
\liminf_{n\rightarrow \infty} \sigma_w(\mathcal{T}_{n,x}) \subset
\liminf\limits_{n\rightarrow \infty} (\bigcup \limits _{x \in
\mathbb{X}}{} \sigma_w(\mathcal{T}_{n,x})). $

\noindent As we have  $\liminf\limits_{n\rightarrow \infty}
\sigma_w(\mathcal{T}_{n,x})= \lim\limits_{n\rightarrow \infty}
\sigma_w(\mathcal{T}_{n,x})= \sigma_w(\mathcal {T}_x), \forall x
\in \mathbb{X},$ then
  $\bigcup \limits _{x \in \mathbb{X}}{} \sigma_w(\mathcal{T}_x)
\subset \liminf\limits_{n\rightarrow \infty}
 (\bigcup \limits _{x \in \mathbb{X}}{}
\sigma_w(\mathcal{T}_{n,x})).$ Since $\bigcup \limits _{x \in
\mathbb{X}}{} \sigma_w(\mathcal{T}_x)= \sigma_w(\mathcal{T})$ and
$\bigcup \limits _{x \in \mathbb{X}}{}
\sigma_w(\mathcal{T}_{n,x})= \sigma_w(\mathcal{T}_n),$ then
  $\sigma_w(\mathcal{T})
\subset \liminf\limits_{n\rightarrow \infty}
\sigma_w(\mathcal{T}_n).$ As we have already from Theorem
\ref{uppersemi}, $\limsup\limits_{n\rightarrow \infty} \sigma_w(
\mathcal{T}_n) \subset \sigma_w(\mathcal{T}),$ then
$\lim\limits_{n\rightarrow \infty} \sigma_w( \mathcal{T}_n) =
\sigma_w(\mathcal{T}).$

\ep

{\small
\noindent Mohammed Berkani,\\
 \noindent Science faculty of Oujda,\\
\noindent University Mohammed I,\\
\noindent Laboratory LAGA, \\
\noindent Morocco\\
\noindent berkanimo@aim.com\\

\end{document}